\input amstex
\magnification=\magstep1 
\baselineskip=13pt
\documentstyle{amsppt}
\vsize=8.7truein \CenteredTagsOnSplits \NoRunningHeads
 \def\conv{\operatorname{conv}}
 \def\PP{\bold{Pr}\thinspace}

 \topmatter

\title  A bound for the number of vertices of a polytope with applications \endtitle
\author Alexander Barvinok  \endauthor
\address Department of Mathematics, University of Michigan, Ann Arbor,
MI 48109-1043, USA \endaddress
\email barvinok$\@$umich.edu \endemail
\date April 2012 \enddate
\thanks  This research was partially supported by NSF Grant DMS 0856640,
by a United States - Israel BSF grant 2006377 and by ERC Advanced Grant No 267195.
\endthanks 
\abstract We prove that the number of vertices of a polytope of a particular kind is exponentially large in the dimension
of the polytope. As a corollary, we prove that an $n$-dimensional centrally symmetric polytope with $O(n)$ facets has 
$2^{\Omega(n)}$ vertices and that the number of $r$-factors in a $k$-regular graph is exponentially large in the number of 
vertices of the graph provided $k \geq 2r+1$ and every cut in the graph with at least two vertices on each side has 
more than $k/r$ edges.
\endabstract
\endtopmatter

\document

\head 1. Introduction and main results \endhead

Let ${\Bbb R}^n$ be Euclidean space with the standard scalar product $\langle \cdot, \cdot \rangle$ and the 
associated Euclidean norm $\| \cdot \|$. 
A polytope $P \subset {\Bbb R}^n$ is the convex hull of a finite set of points. We say that $P$ is $n$-{\it dimensional}
if $P$ has a non-empty interior. The intersection of $P$ with a supporting affine hyperplane is called a {\it face} of $P$.
Faces of $P$ of dimension 0 are called
{\it vertices} and faces of codimension 1 are called {\it facets} of $P$.

In this paper we prove the following result.
\proclaim{(1.1) Theorem} For every $\alpha, \beta \geq 1$  there is $\gamma=\gamma\left(\alpha, \beta \right) >0$ independent of 
$n$ such that the following holds.

Suppose that $P \subset {\Bbb R}^n$ is a polytope containing the set 
$$\Bigl\{x\in {\Bbb R}^n: \quad \left| \langle x, u_i \rangle \right| \ \leq \ 1 \quad \text{for} \quad i=1, \ldots, m\Bigr\}$$
where $\|u_i\| \leq 1$ for $i=1, \ldots, m$ and $m \leq \alpha n$. Suppose further that $P$ lies inside the ball
$$\Bigl\{x \in {\Bbb R}^n: \quad \|x\| \ \leq \ \beta \sqrt{n} \Bigr\}.$$
Then $P$ has at least $\exp\{\gamma n\}$ vertices.
\endproclaim

We obtain the following estimate of $\gamma=\gamma(\alpha, \beta)$ from our proof. 
 Let us choose any $0 < \epsilon < 1$ and let $\rho > 0$ be any number such that the following two inequalities 
 hold:
$$\split &\alpha \ln \left(1-\exp\left\{-{\rho^2 \over 2}\right\} \right) \ > \ -{\epsilon^2 \over 4} \qquad \text{and} \\
&\gamma={(1-\epsilon)^2 \over 2 \beta^2 \rho^2} \left(1-\exp\left\{-{\rho^2 \over 2}\right\} \right) + \alpha  \ln \left(1-\exp\left\{-{\rho^2 \over 2}
\right\}\right) \ > \  0. \endsplit$$
Then $\gamma$ defined by the second formula
satisfies the conclusion of Theorem 1.1 for all sufficiently large $n$.

Our first corollary is a lower bound for the number of vertices of a {\it centrally symmetric} polytope $P$, that is, a 
polytope $P$ satisfying $P=-P$.
\proclaim{(1.2) Corollary} For every $\alpha \geq 1$ there exists $\gamma=\gamma(\alpha) >0$ such that if 
$P$ is an $n$-dimensional centrally symmetric polytope with not more than $\alpha n$ facets then $P$ has at least 
$2^{\gamma n}$ vertices.
\endproclaim

By duality, an $n$-dimensional centrally symmetric polytope with $O(n)$ vertices has $2^{\Omega(n)}$ facets.
In \cite{F+77}  Figiel, Lindenstrauss and Milman, as a corollary of a powerful new method, proved
that for an $n$-dimensional centrally symmetric polytope with
$v$ vertices and $f$ facets the inequality
$$\left(\log v\right) \cdot \left(\log f \right) \ \geq \ \gamma n \tag1.2.1$$
holds for some absolute constant $\gamma >0$. More generally, they obtained (1.2.1) for any polytope, symmetric 
or not, which contains the unit ball and is contained in a ball of radius $O\left(\sqrt{n}\right)$.
We note that if $f=O(n)$ then inequality (1.2.1) implies that 
$v=2^{\Omega(n/\log n)}$ and hence the estimate of Corollary 1.2 is sharper than (1.2.1) in this case. 

Our second application is combinatorial.

Let $G$ be a $k$-regular graph with a finite set $V$ of vertices and the set $E$ of edges. 
Thus every vertex $v \in V$ is incident to precisely 
$k$ edges of $G$ (we do not allow multiple edges or loops). An $r$-regular subgraph $H$ of $G$ with the same set 
$V$ of vertices is called an {\it $r$-factor} of $G$. In particular, a $1$-factor is also known as a {\it perfect matching} in 
$G$. For a set $U \subset V$ of vertices, we denote by $\delta(U) \subset E$ the {\it cut} associated with $U$, that is, the 
set of edges of $G$ with exactly one endpoint in $U$. We denote by $|X|$ the cardinality of a finite set $X$.

We prove that the number of $r$-factors in a $k$-regular graph 
without cuts of small size is exponentially large in the number of vertices of the graph.

\proclaim{(1.3) Corollary} For any positive integers $k$ and $r$ with $k \geq 2r+1$ there is a constant 
$\gamma=\gamma(k, r) >0$ such that the following holds.

Every $k$-regular graph $G$ with vertex set $V$, such that all cuts $\delta(U)$ with
 $2 \leq |U| \leq |V|-2$ have size at least
 $$\bigl| \delta(U)\bigr| > {k \over r},$$
has at least $2^{\gamma |V|}$ distinct $r$-factors.
\endproclaim

Note that the complement to an $r$-factor is a $(k-r)$-factor, so our result also produces an estimate for the number of 
factors of degree greater than one half of the degree of the graph. 

The most tantalizing situation is that of $k=3$ and $r=1$, when Corollary 1.3 asserts that the number of perfect matchings 
of a 3-regular (also known as {\it cubic}) graph is exponentially large in the number $|V|$ of vertices of the graph, provided 
$\big|\delta(U)\big| \geq 4$ as long as $2 \leq |U| \leq |V|-2$. This falls short of the recent result of \cite{E+11}, where it is 
proven that it suffices to have $\big|\delta(U)\big| \geq 2$, and hereby the Lov\'asz-Plummer conjecture is confirmed. 
We hope, however, 
that our method can be sharpened to provide an alternative (and, perhaps, simpler) proof of the conjecture.

We prove Theorem 1.1 and Corollary 1.2 in Section 2 and Corollary 1.3 in Section 3.

 The idea of the proof of Theorem 1.1 is, roughly, as follows. We consider
the maximum of a random linear function on $P$. We argue that if the number of vertices of $P$ is small,
then the maximum is also small. We then argue that if we go from the origin along a random direction
then we stay long enough
inside $P$. This proves that the maximum of a random linear
function on $P$ is large enough and hence $P$ has sufficiently many vertices.
A similar argument is used in Section VI.8 of \cite{Ba02} in the proof of the Figiel-Lindenstrauss-Milman inequality (1.2.1).
M. Rudelson pointed out to the author that a proof of Theorem 1.1, although with a weaker bound for $\gamma$, can also be obtained
using volume estimates of \cite{Gl88}. 

To prove Corollary 1.2, we apply a linear transformation so that the image of $P$ satisfies the conditions of Theorem 1.1.

To prove Corollary 1.3,  we consider a polytope $P_r(G)$ whose vertices correspond to 
$r$-factors of $G$ and then apply Theorem 1.1. N. Alon pointed out to the author that if $k$ and $r$ are both even
then, using an Eulerian trail in $G$, one can deduce the existence of exponentially many distinct $r$-factors in a $k$-regular 
graph from the Van der Waerden  bound in the bipartite case (see, for example, Theorem 8.1.3 of \cite{LP86}), just as one can deduce the existence of a single $r$-factor in Petersen's Theorem, see Theorem 6.2.4 in \cite{LP86}. An attempt 
to apply a similar decomposition argument to other values of $k$ and $r$ seems to require a higher connectivity of $G$ 
than that assumed by Corollary 1.3. For example, for $k=5$ and $r=2$ no combinatorial argument seems to be readily available.

The paper \cite{BS07} describes a general method of asymptotic counting of combinatorial structures through optimization 
of a random linear function.

\head 2. Proofs of Theorem 1.1 and of Corollary 1.2 \endhead 

Let us fix the standard Gaussian probability measure in ${\Bbb R}^n$ with  density 
$${1 \over (2 \pi)^{n/2}} \exp\left\{-{\|x\|^2 \over 2}\right\} \quad \text{for} \quad x \in {\Bbb R}^n.$$
\proclaim{(2.1) Lemma} 
\roster 
\item For any $0 < \epsilon < 1$ we have 
$$\PP\Bigl(y \in {\Bbb R}^n:\   \ \|y\|^2 \ \leq \ (1-\epsilon) n \Bigr) \ \leq \ \exp\left\{ - {\epsilon^2 n \over 4} \right\}.$$
\item 
Let $a \in {\Bbb R}^n$ be a point, then for any $\tau \geq 0$
$$\PP\Bigl(y \in {\Bbb R}^n: \ \big\langle y, a\big\rangle\ \geq \ \tau \Bigr) \ \leq \ 
{1 \over 2} \exp\left\{ -{\tau^2 \over 2 \|a\|^2} \right\}.$$
\item For any $\rho \geq 0$ and any vectors $u_1, \ldots, u_m \in {\Bbb R}^n$ such that 
$\|u_i\| \leq 1$ for $i=1, \ldots, m$, we have 
$$\PP\Bigl( y \in {\Bbb R}^n: \ \left|\langle u_i, y \rangle \right| \leq \rho \quad \text{for} \quad i=1, \ldots, m\Bigr) \ \geq \ 
\left(1- \exp\left\{-{\rho^2 \over 2}\right\}\right)^m.$$
\endroster
\endproclaim
\demo{Proof} 
The inequality of Part (1) can be found, for example, in Corollary V.5.5 of \cite{Ba02}.

The function $y \longmapsto \langle y, a \rangle$ is a centered Gaussian random 
variable with variance $\|a\|^2$ and Part (2) follows by the standard Gaussian tail estimate.

By the \v Sidak Lemma, see, for example, \cite{Ba01}, we have 
$$\PP\Bigl(y \in {\Bbb R}^n: \ \left| \langle u_i, y \rangle \right| \leq \rho \quad \text{for} \quad i=1, \ldots, m \Bigr) \ \geq \ 
\prod_{i=1}^m \PP \Bigl(y \in {\Bbb R}^n: \ \left| \langle u_i, y \rangle \right| \leq \rho \Bigr)$$
(informally, the \v Sidak Lemma states that ``slabs" are positively correlated with respect to the Gaussian measure).
Since $y \longmapsto \langle u_i, y \rangle$ is a centered Gaussian random variable of variance $\|u_i\|^2 \leq 1$, 
the proof of Part (3) follows from Part (2).
{\hfill \hfill \hfill} \qed
\enddemo

\subhead (2.2) Proof of Theorem 1.1 \endsubhead 

Let us choose any $0 < \epsilon < 1$ and a sufficiently large $\rho=\rho(\alpha, \beta, \epsilon) >0$ such that the following two inequalities hold:
$$\alpha \ln \left(1 - \exp\left\{-{\rho^2 \over 2} \right\} \right) \ > \ -{\epsilon^2 \over 4} \tag2.2.1$$
and
$$\gamma={(1-\epsilon)^2  \over 2\beta^2 \rho^2} \left(1 - \exp\left\{ -{\rho^2 \over 2} \right\} \right) + \alpha \ln \left(1-\exp\left\{-{\rho^2 \over 2}\right\}
\right)> 0 \tag2.2.2$$
We prove that the conclusion of Theorem 1.1 holds with $\gamma$ defined by (2.2.2) for all sufficiently large 
$n > n_0(\alpha, \beta, \epsilon, \rho)$.

Let us consider the polyhedron
$$Q=\Bigl\{y \in {\Bbb R}^n: \ \left| \langle y, u_i \rangle \right| \ \leq \ \rho \quad 
\text{for} \quad i=1, \ldots, m \Bigr\}.$$
By Part (3) of Lemma 2.1 we have 
$$\PP\Bigl(y: \ y \in Q \Bigr) \ \geq \ \left(1 - \exp\left\{ -{\rho^2 \over 2}\right\} \right)^m
\ \geq \ \left(1 - \exp\left\{ -{\rho^2 \over 2}\right\} \right)^{\alpha n}. $$
By Part (1) of Lemma 2.1 and (2.2.1) we conclude that 
$$\PP\Bigl(y: \ y \in Q \quad \text{and} \quad \|y\|^2 \ \geq \ (1-\epsilon)n \Bigr) \ \geq \ 
{1 \over 2} \left(1 - \exp\left\{ -{\rho^2 \over 2}\right\} \right)^{\alpha n} \tag2.2.3$$
for all sufficiently large $n> n_0(\alpha, \beta, \epsilon, \rho)$.

We consider the maximum value of the linear function $x \longmapsto \langle x, y \rangle$ on $P$.
Since for every $y \in Q$ we have $\rho^{-1} y \in P$ we conclude that 
$$\max_{x \in P} \langle x, y \rangle \ \geq \ \left\langle \rho^{-1} y, \ y \right\rangle = \rho^{-1}  \|y\|^2  
\quad \text{for all} \quad y \in Q. $$
Therefore, from (2.2.3) we have
$$\PP\left(y: \ \max_{x \in P} \langle x, y \rangle \ \geq \ \rho^{-1} (1-\epsilon) n \right) \ \geq \ {1 \over 2}
\left(1 -\exp\left\{ -{\rho^2 \over 2}\right\}\right)^{\alpha n} \tag2.2.4$$
for all sufficiently large $n> n_0(\alpha, \beta, \epsilon, \rho)$.

On the other hand, the maximum value of a linear function on a polytope is attained, in particular, at a vertex 
of $P$. Therefore, taking $W$ to be the set of vertices of $P$, from Part (2) of Lemma 2.1, we conclude that 
$$\split \PP\left(y: \ \max_{x \in P} \langle x, y \rangle \ \geq \ \tau \right) \ \leq \ &\sum_{a \in W}
\PP\bigl(y: \ \langle y, a \rangle \ \geq \ \tau \bigr) \ \leq \ {1 \over 2} \sum_{a \in W} \exp\left\{- {\tau^2 \over 2\|a\|^2}\right\}
\\ \leq \ &{|W| \over 2} \exp\left\{ -{\tau^2 \over 2\beta^2 n}\right\}
.\endsplit$$
Substituting 
$$\tau=\rho^{-1} (1-\epsilon) n,$$
we obtain
$$ \PP\left(y: \ \max_{x \in P} \langle x, y \rangle \ \geq \ \rho^{-1}(1-\epsilon) n \right) 
\ \leq \ {|W| \over 2} \exp\left\{ -{(1-\epsilon)^2 n \over 2 \rho^2  \beta^2}\right\}. \tag2.2.5$$
Combining (2.2.5) and (2.2.4) and using (2.2.2), we conclude that 
$$|W| \ \geq \ \exp\left\{ {(1-\epsilon)^2 n \over 2 \rho^2 \beta^2}\right\} \left( 1- \exp\left\{ -{\rho^2 \over 2}\right\} \right)^{\alpha n} 
\ \geq \ \exp\left\{ \gamma n \right\}$$
for all sufficiently large $n$, as desired.
{\hfill \hfill \hfill} \qed

\subhead (2.3) Proof of Corollary 1.2 \endsubhead

We can write 
$$P=\Bigl\{x \in {\Bbb R}^n: \quad \left|\langle u_i, x \rangle\right|\ \leq \ \delta_i  \quad \text{for} \quad i=1, \ldots, m\Bigr\},$$
where $u_1, \ldots, u_m$ are the unit normals to the facets of $P$ and $\delta_i >0$.
Applying to $P$ an invertible linear transformation, we 
may assume, additionally, that $P$ contains the unit ball and is contained in the ball of radius $\sqrt{n}$,
where both balls are centered at the origin (see, for example, Sections V.2 and VI.8 of \cite{Ba02}). Since $P$ contains the unit ball,
we must have $\delta_i \geq 1$ for all $i=1, \ldots, m$ and the proof follows by Theorem 1.1.
{\hfill \hfill \hfill} \qed

\head 3. Proof of Corollary 1.3 \endhead

\subhead (3.1) The polytope \endsubhead Let $G$ be a graph with the set $V$ of vertices and the set $E$ of edges.
We denote by ${\Bbb R}^E$ the Euclidean space of all real-valued functions $x: E \longrightarrow {\Bbb R}$. 
We use the standard scalar product 
$$\langle x, y \rangle = \sum_{e \in E} x(e) y(e) \quad \text{for all} \quad x, y \in {\Bbb R}^E$$ and 
the corresponding Euclidean norm $\|x\| = \sqrt{\langle x, x \rangle}$.

For a subset $H\subset E$ we consider a vector (indicator function) $[H] \in {\Bbb R}^E$ defined as follows:
$$[H](e)=\cases 1 &\text{if $e$ is an edge of $H$} \\ 0 &\text{otherwise.} \endcases$$
We define the {\it $r$-factor polytope} $P_r(G)$ as the convex hull
$$P_r(G) =\conv\Bigl([H]: \quad H \text{\ is an $r$-factor of\ } G \Bigr).$$

We will need the following description of $P_r(G)$ by a system of linear inequalities (3.1.1)--(3.1.3), see Corollary 33.2a of 
\cite{Sc03} :
$$0 \ \leq \ x(e) \ \leq \ 1 \quad \text{for all} \quad e \in E, \tag3.1.1$$
$$\sum_{e \in \delta(v)} x(e) =r \quad \text{for all} \quad v \in V,  \tag3.1.2$$ and 
$$\aligned \sum_{e \in \delta(U) \setminus F} x(e) - \sum_{e \in F} x(e) \ \geq \ &1- |F| \quad \text{for all} \quad U \subset V,\ F \subset \delta(U) \\ &  \text{such that} \quad r|U| +|F| \quad \text{is odd}. \endaligned \tag3.1.3$$
Our first goal is to show that if $G$ is $k$-regular graph without small cuts then the vector $a \in {\Bbb R}^E$ 
with $a(e)=r/k$ for all $e \in E$ lies sufficiently deep inside the polytope $P_r(G)$.
\proclaim{(3.2) Lemma} Suppose that $G$ is $k$-regular, that $k \geq 2r+1$ and that 
$$\big| \delta(U) \big| \  > \ {k \over r}$$
for all $U \subset V$ such that $2 \leq |U| \leq |V|-2$.
Let us define $\epsilon=\epsilon(k, r) >0$ by
$$\epsilon =\min\left\{{r \over k} -{1 \over \lceil {k+1 \over r}\rceil }, \quad {1 \over 2k} \right\}.$$
 Let $a \in {\Bbb R}^E$ be the vector such that 
$$a(e)={r \over k} \quad \text{for all} \quad e \in E$$
and let $y \in {\Bbb R}^E$ be a vector such that 
$$\sum_{e \in \delta(v)} y(e)=0 \quad \text{for all} \quad v \in V$$
and 
$$\left| y(e) \right| \ \leq \ \epsilon \quad \text{for all} \quad e \in E.$$
Then for $x =a+y$ we have $x \in P_r(G)$.
\endproclaim
\demo{Proof} Clearly, vector $x$ satisfies (3.1.1) and (3.1.2). 
Moreover,
$$ {2r-1 \over 2k} \ \leq \ x(e) \ \leq \ {2r+1 \over 2k} \quad \text{for all} \quad e \in E.$$
If in (3.1.3) we increase $|F|$ by 1 then the left hand side decreases at most by $(2r+1)/k$ while the right hand side 
decreases by 1. 
Therefore, it suffices to check (3.1.3) when $|F|=0$.
 Furthermore, if $|U|=1$ or if $|V\setminus U|=1$, inequality 
(3.1.3) follows by (3.1.2).

If $|F|=0$ then the left hand side of (3.1.3) is at least 
$$\big|\delta(U)\big| \Bigl({r \over k}-\epsilon\Bigr) \ \geq \ 1$$ 
and (3.1.3) holds.
{\hfill \hfill \hfill} \qed
\enddemo
 
\subhead (3.3) Proof of Corollary 1.3 \endsubhead 

All implied constants in the $O(\cdot)$ and $\Omega(\cdot)$ terms below may depend on $k$ and $r$, but not on 
$n=|V|$.

Since $G$ is $k$-regular, we have $|E|=k|V|/2$. Let $L \subset {\Bbb R}^E$ be the subspace defined 
by the equations
$$\sum_{e \in \delta(v)} x(e) =0 \quad \text{for all} \quad v \in V.$$
Hence 
$$n=\dim L \ \geq \ |E| -|V| = \left({k \over 2} -1\right) |V| = \Omega(V).$$
We identify $L$ with ${\Bbb R}^n$. Let $P=P_r(G) -a$, where $a$ is the vector of Lemma 3.2. 
Then $P \subset {\Bbb R}^n$ by (3.1.2). Since 
$$\| [H]\| \ = \sqrt{ r|V| \over 2}$$
for any $r$-factor $H$ of $G$ and 
$$\|a\| ={r \over k} \sqrt{k |V| \over 2},$$
we conclude that $P$ lies in a ball of radius $O\left(\sqrt{n}\right)$ centered at the origin.
 
Moreover, by Lemma 3.2, the polytope $P$ contains the set 
$$\Bigl\{ x \in {\Bbb R}^n: \quad \left| \langle u_e, x \rangle \right| \ \leq \ \epsilon \quad \text{for all} \quad e \in E \Bigr\},$$
where $u_e$ is the orthogonal projection of $[e]$ onto $L$. In particular, $\|u_e\| \leq 1$ for all $e \in E$.
Since $|E| =O(n)$ and $\epsilon =\Omega(1)$, the proof is obtained by applying Theorem 1.1 to the dilated polytope
$\epsilon^{-1} P$.
{\hfill \hfill \hfill} \qed

\head Acknowledgments \endhead

The author is grateful to Alex Samorodnitsky for many useful discussions and comments on a draft 
of this paper, to Imre B\'ar\'any and Alfr\'ed R\'enyi Institute of Mathematics (Budapest) for hospitality
during his work on this paper, to Mark Rudelson for pointing out to connections with \cite{Gl88}, to Noga Alon 
for explaining graph decomposition techniques and to anonymous referees for careful reading of the 
paper and helpful suggestions.

\enddocument
\end